\documentclass{eajam}
\renewcommand\thisnumber{x}
\renewcommand\thisyear {200x}
\renewcommand\thismonth{xxx}
\renewcommand\thisvolume{xx}

\usepackage{charter}
\usepackage[charter]{mathdesign}
\usepackage{pdfpages,color,tikz,pdfsync,enumitem}
\usepackage{graphicx}

\allowdisplaybreaks[4]

\numberwithin{equation}{section}

\newcommand{\R}{\mathbb{R}}
\newcommand{\J}{\mathbb{J}}
\newcommand{\mC}{\mathbb{C}}
\newcommand{\Z}{\mathbb{Z}}

\newcommand{\T}{\mathbb{T}^2}

\newcommand{\re}{\mathrm{Re}\,}
\newcommand{\im}{\mathrm{Im}\,}

\newcommand{\Hess}{\mathrm{Hess}}

\newcommand{\bvec}[1]{\boldsymbol{#1}}

\newcommand{\vj}{\bvec{j}}

\newcommand{\vq}{\bvec{q}}

\newcommand{\vx}{\bvec{x}}
\newcommand{\va}{\bvec{a}}

\newcommand{\val}{\bvec{a}^0}
\newcommand{\vb}{\bvec{b}}

\newcommand{\p}{\partial}

\newcommand{\Div}{\nabla\cdot}

\newcommand{\la}{\left\langle}
\newcommand{\ra}{\right\rangle}
\newcommand{\nn}{\nonumber\\}

\newcommand{\SE}{Schr\"{o}dinger equation}

\newcommand{\wto}{\rightharpoonup}

\newcommand{\be}{\begin{equation}}
\newcommand{\ee}{\end{equation}}

\newcommand{\M}{W^{-1,1}(\T)}
\newcommand{\eprint}{arxiv:}

\begin{document}

\markboth{Y.~Zhu}{Quantized Vortex Dynamics of  NLSW on the Torus}
\title{Quantized Vortex Dynamics of  the  Nonlinear Schr\"{o}dinger Equation with Wave Operator on the Torus}

\author[Zhu]{Yongxing Zhu\corrauth}
\address{Department of Mathematical Sciences, Tsinghua University, Beijing, 100084, China}
\email{{\tt zhu-yx18@mails.tsinghua.edu.cn} (Yongxing Zhu)}

\begin{abstract}
We derive rigorously the reduced dynamical law for quantized vortex dynamics of the nonlinear \SE{} with wave operator on the torus when the core size of vortex $\varepsilon\to 0$.  It is proved that  the reduced dynamical law of the nonlinear \SE{} with wave operator is a mixed state of the vortex motion laws for the nonlinear wave equation and  the nonlinear \SE. We will also investigate the convergence of the reduced dynamical law of the nonlinear \SE{} with wave operator to the vortex motion law of the nonlinear \SE{} via numerical simulation.
\end{abstract}



\maketitle

\section{Introduction}\label{sec:Introduction}

In this paper, we study the vortex dynamics of the nonlinear \SE{} with wave operator (NLSW) \cite{BaoCai2014NSW,ColinFabrie1998NSW,BergeColin1995NSWinplasma,Zhao2017MultiscaleTimeIntegratorNSW}:
\begin{equation}\label{eq:NSW}
 -i\p_t u^\varepsilon+\mu k_\varepsilon\p^2_tu^\varepsilon-\Delta u^\varepsilon+\frac{1}{\varepsilon^2}(|u^\varepsilon|^2-1)u^\varepsilon=0,\vx\in\T,t>0,
\end{equation}
with initial data
\begin{equation}
u^\varepsilon(\vx,0)=u^\varepsilon_0(\vx),\quad \p_tu^\varepsilon(\vx,0)=u^\varepsilon_1(\vx),\quad \vx\in\T.
\end{equation}
Here, $0<\varepsilon<<1$ is a parameter related to the vortex core size, $k_\varepsilon=1/|\log\varepsilon|$, $0<\mu<<1$ is a positive parameter, $\T=(\R/\Z)^2$ is the unit torus, $\vx=(x,y)^T$ is the spatial coordinate, $t$ denotes the time variable, $u^\varepsilon$ is a complex-valued function which is called the order parameter, and $u^\varepsilon_0,u^\varepsilon_1$ are  complex-valued initial data.

We define the Ginzburg-Landau functional (energy) by \cite{Lin1999VortexDynamicsWave,Jerrard1999GinzburgLandauwave,Yu2011VortexDynamicsKleinGordon}:
\begin{equation}\label{eq:def of E}
E^\varepsilon(u^\varepsilon(t)):=\int_{\T}e^\varepsilon(u^\varepsilon(\vx,t))d\vx,
\end{equation}
the momentum by \cite{Jerrard1999GinzburgLandauwave,Yu2011VortexDynamicsKleinGordon}:
\begin{equation}\label{eq:def of Q}
  \bvec{Q}(u^\varepsilon(t)):=\int_{\T}\vj(u^\varepsilon(\vx,t))d\vx,
\end{equation}
where the energy density $e^\varepsilon(v)$, the current $\vj(v)$ and the Jacobian $J(v)$ are defined as follows: for any complex-valued function $v:\T\to\mC^2$,
\begin{equation}\label{eq:def of e j J}
\begin{aligned}
&\vj(v):=\im (\overline{v}\nabla v), \quad
e^\varepsilon(v):=\frac{1}{2}|\nabla v|^2+\frac{1}{4\varepsilon^2}(1-|v|^2)^2,\\ &J(v)=\frac{1}{2}\nabla\cdot(\J \vj(v))=\im (\partial_x\overline{v}\,\partial_y v),
\end{aligned}
\end{equation}
with $\overline{v},\im v$ denoting the complex conjugate and imaginary part of $v$, respectively, and 
\begin{equation}
  \J=\left(\begin{array}{cc}0&1\\-1&0 \end{array}\right).
\end{equation}

NLSW \eqref{eq:NSW} arises in the study of the nonrelativistic limit of the nonlinear Klein-Gordon equation \cite{MachiharaNakanishiOzawa2002LimitOfNKG,Schoene1979NonrelativisticLimitKGandD,Masayoshi1984NonrelativisticApproximationNKG}, the Langmuir wave envelope approximation in plasma \cite{ColinFabrie1998NSW,BergeColin1995NSWinplasma}, and  the modulated planar pulse approximation of the sine-Gordon equation for light bullets \cite{BaoDongXin2010sGperturbedNS,Xin2000ModellingLightBullets}. 
It was shown that as $\mu\to 0$, NLSW \eqref{eq:NSW} converges to the nonlinear \SE{} (NLS) \cite{BaoCai2014NSW,BergeColin1995NSWinplasma,MachiharaNakanishiOzawa2002LimitOfNKG,Schoene1979NonrelativisticLimitKGandD}, which is widely used as a phenomenological model for superfluidity and Bose-Einstein condensate \cite{BaoDuZhang2006RotatingBEC,LinXin1999,Anderson2010ExperimentsVorticeSuperfluid}:
\begin{align}\label{eq:SE}
-i\p_tu^\varepsilon-\Delta u^\varepsilon+\frac{1}{\varepsilon^2}(|u^\varepsilon|^2-1)u^\varepsilon=0,\quad \vx\in\T,t>0.
\end{align}
A key feature of superfluidity is the appearance of quantized vortex \cite{Anderson2010ExperimentsVorticeSuperfluid}. Mathematically, quantized vortex means the topological defect of the order parameter, which is the zero point with nonzero winding number (degree) in two-dimensional cases. Quantized vortex arises in many physical phenomena, such as superfluidity, superconductivity and Bose-Einstein condensate, and is widely observed in experiments \cite{AransonKramer2002ComplexGinzburgLandau,Anderson2010ExperimentsVorticeSuperfluid}. From both analytical and numerical study \cite{Mironescu1995StabilityRadialSolutionGinzburgLandau,Anderson2010ExperimentsVorticeSuperfluid,BaoZengZhang2008WaveVortex}, a vortex is stable only when its degree is $1$ or $-1$. 

 The interaction and dynamics of quantized vortex were widely studied in the past decades. We refer to \cite{LinXin1999,test,Rubinstein1991SelfinducedmotionLinedefects,Neu1990VorticesComplexScalarfield,E1994DynamicsGinzburgLandau,BaoTang2013NumericalQuantizedVortexGinzburgLandau,BaoTang2014NumericalStudyQuantizedVorticesSchrodinger,SandierSerfaty2004EstimateGinzburgLandau,Lin1996DynamicGinzburgLandau,BetheulOrlandiSmets2006GinzburgLandautoMotionByMeanCurvature,ChenSternberg2014VortexDynamicsManifold,DuJuGInzburgLandauSphereNumerical,JianLiu2006GinzburgLandauVortexMeanCurvatureFlow,Jerrard1999GinzburgLandauwave,JerrardSoner1998DynamicsGinzburgLandauVortices,JerrardSpirn2008RefinedEstimateGrossPitaevskiiVortices,SandierSerfaty2004GammaConvergenceGinzburgLandau,LinXin1999DynamicGinzburgLandauPlane,Lin1999VortexDynamicsWave,ZhangBaoDu2007SimulationVortexGinzburgLandauSchrodinger,BetheulJerrardSmets2008NLSVortexdynamicsPlane,BaoDuZhang2006RotatingBEC} and references therein for the results on the whole plane $\R^2$ or on a bounded domain with Neumman or Dirichlet boundary condition. In particular, the vortex dynamics of \eqref{eq:NSW} on a bounded domain was studied in \cite{Yu2011VortexDynamicsKleinGordon,Yu2011VortexDynamicsMaxwellKG}. Yu \cite{Yu2011VortexDynamicsKleinGordon} proved that only when $\mu$ is a fixed positive number, we can obtain the reduced dynamical law different from those of the nonlinear wave equation and NLS \eqref{eq:SE}. It was also shown in \cite{Yu2011VortexDynamicsKleinGordon} that for \eqref{eq:NSW} on a bounded domain $\Omega$  with proper $u_1^\varepsilon$ and boundary condition, if the initial data $u^\varepsilon_0$ satisfies
\begin{align}
J(u^\varepsilon_0)\to \pi\sum_{j=1}^{M}d_j\delta_{\va_j^0}(\vx),\quad \text{in}\ W^{-1,1}(\Omega):=[C_0^1(\Omega)]',\  \text{as}\ \varepsilon\to 0,
\end{align}
where $\va_1^0,\cdots,\va_M^0\in\Omega$ are distinct points, $d_1,\cdots,d_M\in\{\pm 1\}$ and $\delta$ is the Dirac-$\delta$ function while $\delta_{\vx_0}(\vx):=\delta(\vx-\vx_0)$.
Then there exist a time $T$ and $\va_j(t),1\le j\le M$, such that the solution $u^\varepsilon(t)$ also satisfies for any $t<T$,
\begin{equation}
J(u^\varepsilon(t))\to \pi\sum_{j=1}^{M}d_j\delta_{\va_j(t)}(\vx),\quad \text{in}\ W^{-1,1}(\Omega),\  \text{as}\ \varepsilon\to 0,
\end{equation}
Moreover, $\va=\va(t)=(\va_1(t),\cdots,\va_M(t))^T$ satisfies
\begin{equation*}
	\mu\ddot{\va}+d_j\J\dot{\va}_j=-\frac{1}{\pi}\nabla W_\Omega(\va),
\end{equation*}
where $W_\Omega$ is the renormalized energy on $\Omega$ of the following form:
\begin{equation}\label{eq:RDL on Omega}
  W_\Omega(\va)=-\pi\sum_{1\le j\ne k\le M}d_jd_k\log|\va_j-\va_k|+\text{Reminding term determined by} \ \Omega.
\end{equation}

For \eqref{eq:NSW} on the torus, there are two differences: (i) The summation of degrees of vortices must be zero since $\T$ is a compact manifold \cite{CollianderJerrard1999GLvorticesSE}. Without loss of generality, we can assume $M=2N $ for integer $N$ and
\begin{equation}
d_1=\cdots=d_N=1,\quad d_{N+1}=\cdots=d_{2N}=-1.
\end{equation}
(ii) The vortex dynamics on the torus is determined by both the positions of vortices and the limit momentum $\bvec{Q}_0:=\lim_{\varepsilon\to 0}\bvec{Q}(u^\varepsilon_0)$ \cite{zhubaojian2023quantized,zhu2023quantized}. And if $u^\varepsilon$ satisfies
\begin{equation}\label{con:convergence of J varphi}
 J (u_0^\varepsilon)\to \pi\sum_{j=1}^{2N} d_j \delta_{\val_j} \ \text{in}\ \M=[C^1(\T)]', \text{as}\  \varepsilon\to 0^+,
 \end{equation}
 then \cite{zhu2023quantized}
 \begin{equation}
 	\bvec{Q}_0\in 2\pi\J\sum_{j=1}^{2N}\va_j^0+2\pi\Z^2.
 \end{equation}

The main purpose of this article is to extend the results in \cite{Yu2011VortexDynamicsKleinGordon} and obtain the reduced dynamical law of \eqref{eq:NSW} involving the influence of momentum. 

Before the statement of our result, we have to introduce some notations. We define 
\begin{equation*}
	(\T)_*^{2N}=\{(\vx_1,\cdots,\vx_{2N})^T\in(\T)^{2N}|\vx_k\ne\vx_l \ \text{for}\ 1\le k<l\le 2N\}.
\end{equation*}
For any $\va=(\va_1,\cdots,\va_{2N})^T\in(\T)_*^{2N}$ and $\vq\in 2\pi\sum_{j=1}^{2N}d_j\va_j+2\pi\Z^2$, the renormalized energy is defined by \cite{IgnatJerrard2021Renormalizedenergymanifold}
\begin{equation}\label{eq:def of W}
	W(\va;\vq):=-\pi\sum_{1\le k\ne l\le 2N}d_kd_lF(\va_k-\va_l)+\frac{1}{2}\left|\vq\right|^2,
\end{equation}
where $F$ is the solution of 
\begin{equation}\label{eq:define of F}
	\Delta F(\vx)=2\pi(\delta(\vx)-1),\vx\in\T,\quad\text{with}\  \int_{\T}Fd\vx=0.
\end{equation}
We define
\begin{equation}\label{eq:def of I}
    \gamma:=\lim_{\varepsilon\to 0}\left(\inf_{u\in H^1_g(B_1(\bvec{0}))}\int_{B_1(\bvec{0})}e^\varepsilon(u)d\vx-\pi\log\frac{1}{\varepsilon}\right),
\end{equation}
where $H^1_g(B_1(\bvec{0}))$ is a function space defined by
\begin{equation*}
  H_g^1(B_1(\bvec{0}))=\left\{u\in H^1(B_1(\bvec{0}))\left| u(\vx)=g(\vx)=\frac{x+iy}{|\vx|}\ \text{for}\  \vx\in \p B_1(\bvec{0})\right.\right\}.
\end{equation*}
 Then we define 
\begin{equation}\label{eq:def of We}
  W_\varepsilon(\va;\vq):=2N\left(\pi\log\frac{1}{\varepsilon}+\gamma \right)+W(\va;\vq).
\end{equation}
Our main result is stated as follows:

\begin{theorem}[Reduced dynamical law of NLSW]\label{thm:dynamics NSW}
Assume there exist $\va^0=(\va_1^0,\dots,$ $\va^0_{2N})^T\in(\T)^{2N}_*$, $\bvec{q}_0\in 2\pi\sum_{j=1}^{2N}d_j\va_j^0+2\pi\Z^2$ such that the initial data of \eqref{eq:NSW} satisfies \eqref{con:convergence of J varphi} and 
 \begin{align}
    &\bvec{Q}_0=\J\vq_0,\quad
    \lim_{\varepsilon\to 0}(E^\varepsilon(u^\varepsilon_0)- W_\varepsilon(\val;\vq_0))=0,\quad k_\varepsilon\int_{\T}|u_1^\varepsilon|^2d\vx\to0.\label{Con:strong bd of energy}
\end{align}
Then there exist  Lipschitz paths $\va_j:[0,T)\to\T$, such that
\begin{equation*}
    J(u^\varepsilon(\vx,t))\to \pi\sum_{j=1}^{2N}d_j\delta_{\va_j(t)},\quad k_\varepsilon e^\varepsilon(u^\varepsilon(\vx,t))\to\pi\sum_{j=1}^{2N}\delta_{\va_j(t)},
\end{equation*}
and $\va=\va(t)=(\va_1(t),\cdots,\va_{2N}(t))^T$ satisfies 
\begin{equation}\label{eq:NSWODE}
\mu\ddot{\va}_j+ d_j\J\dot{\va}_j=-\frac{1}{\pi}\nabla_{\va_j}W(\va;\vq_*(\va)),
\end{equation}
where $\vq_*(\va)=\vq_*(\va(t))$ is continuous and satisfies \begin{equation}\label{eq:def of q*}
  \vq_*(\va(t))\in2\pi\sum_{j=1}^{2N}d_j\va_j(t)+2\pi\Z^2,
  \end{equation}
and the initial data is given by
\begin{equation}\label{eq:initial ODE}
  \va_j(0)=\va_j^0,\quad \dot{\va}_j(0)=0,\quad \vq_*(\va(0))=\vq_0.
\end{equation}
\end{theorem}

{The rest of this  paper is organized as follows: In Section \ref{sec:Preliminaries}, we will introduce the canonical harmonic map and renormalized energy on the torus. 
In Section \ref{sec:NSW}, we will give the proof of Theorem \ref{thm:dynamics NSW}. We first prove the convergence of the Jacobian and the current of $u^\varepsilon$. Then we give some refined lower bounds of the energy of $u^\varepsilon$. Finally, we give the prove of the reduced dynamical law. In Section \ref{sec:numerical}, we will investigate the convergence of the reduced dynamical law of NLSW \eqref{eq:NSW} to the vortex motion law of NLS \eqref{eq:SE} as $\mu\to 0$ via numerical simulation.}

\section{Notations and conservation laws}\label{sec:Preliminaries}

\subsection{Canonical harmonic map and renormalized energy}
Recall that $F$ is the solution of \eqref{eq:define of F}. Then following \cite{CollianderJerrard1999GLvorticesSE,zhubaojian2023quantized}, for any $\va=(\va_1,\cdots,\va_{2N})^T$ $\in(\T)_*^{2N}$ and $\vq\in2\pi\sum_{j=1}^{2N}d_j\va_j+2\pi\Z^2$, we can define the canonical harmonic map $H=H(\vx;\va,\vq)\in C^\infty_{loc}(\T_*(\va))\cap W^{-1,1}(\T)$ by solving
\begin{align}
|H|=1,\quad\Div \vj(H)=0,\quad 2J(H)=\nabla\cdot(\J \vj(H))=2\pi\sum_jd_j\delta_{\va_j},\quad \int_{\T}\vj(H)=\bvec{q}.\label{eq:div of jH}
\end{align}

We use $\Hess(v)$ to denote the Hessian matrix of $v$.
Define 
\begin{equation}\label{eq:def of ra}
r(\va)=\frac{1}{4}\min_{1\le j<k\le 2N}|\va_j-\va_k|.
\end{equation}
Then for any $\rho<r(\va)$, we can define
\begin{equation}
\T_\rho(\va)=\T\setminus\cup_{j=1}^{2N}B_\rho(\va_j)
,\quad
\T_*(\va)=\cup_{\rho>0}\T_\rho(\va)=\T\setminus\{\va_1,\cdots,\va_{2N}\}.
\end{equation}
For two complex vectors $\bvec{z}=(z_1,\cdots,z_M)^T,\bvec{w}=(w_1,\cdots,w_M)^T\in\mC^M,M=1,2$, we define
\begin{equation*}
    \la \bvec{z},\bvec{w}\ra:=\re\sum_{j=1}^M\overline{z}_jw_j.
\end{equation*}
With the above notations, Lemma 2.3 in \cite{zhubaojian2023quantized} gives that
for any $\eta\in C_0^2(B_{r(\va)}(\va))$ which is linear in a neighborhood of $\va_{j}$,
we have
\begin{align}
&\int_{\T} \la\Hess(\eta)\vj(H),\J\vj(H)\ra d \vx
=-\nabla \eta(\va_{j})\cdot(\J\nabla_{\va_{j}}W(\va;\vq)).\label{eq:production of jH and eta GP}
\end{align}

From (2.21) in \cite{zhu2023quantized}, for  $\rho<r(\va)$, $W(\va;\vq)$ is a Lipschitz function with respect to $\va$ in $(\T)^{2N}_\rho$ with
\begin{equation}
(\T)^{2N}_\rho:=\{(\vx_1,\cdots,\vx_{2N})\in(\T)^{2N}||\vx_j-\vx_k|>\rho,\forall 1\le j<k\le 2N\}.
\end{equation}
i.e. there exists a constant $C_\rho$ such that 
\begin{equation}\label{eq:Lip W}
\|W(\va;\vq)\|_{C^1((\T)^{2N}_\rho)}\le C_\rho.
\end{equation}

\subsection{Conservation laws}
In this subsection, we list some equalities related to \eqref{eq:NSW}.

With the definition of Hamiltonian $h^\varepsilon_\mu(u^\varepsilon)$
\cite{Yu2011VortexDynamicsKleinGordon}:
\begin{equation}\label{eq:def of hmu}
  h_\mu^\varepsilon(u^\varepsilon):=\frac{\mu k_\varepsilon}{2}|u^\varepsilon_t|^2+e^\varepsilon(u^\varepsilon),
\end{equation}
we have that the solution $u^\varepsilon$ of Equation \eqref{eq:NSW} satisfies the following equalities \cite{Yu2011VortexDynamicsKleinGordon}:
\begin{align}
\frac{\p}{\p t}h^\varepsilon_\mu(u^\varepsilon(\vx,t))&=\nabla\cdot(\re(\overline{u^\varepsilon_t}\nabla u^\varepsilon)),\label{eq:derivative of h}\\
\nabla\cdot \vj(u^\varepsilon(\vx,t))&=\mu k_\varepsilon\p_t\im(\overline{u^\varepsilon}u^\varepsilon_t)- \frac{1}{2} \p_t|u^\varepsilon|^2.\label{eq:div of ju}
\end{align}
A direct corollary of \eqref{eq:derivative of h} is 
\begin{equation}\label{eq:hamiltonian conservation}
H_\mu^\varepsilon(u^\varepsilon(t)):=\int_{\T}\frac{\mu k_\varepsilon}{2}|u^\varepsilon_t(\vx,t)|^2d\vx+E^\varepsilon(u^\varepsilon(t))\equiv \int_{\T}\frac{\mu k_\varepsilon}{2}|u^\varepsilon_1|^2d\vx+E^\varepsilon(u_0^\varepsilon).
\end{equation}
For any $\varphi\in C^\infty(\T)$, we have \cite{Yu2011VortexDynamicsKleinGordon}
\begin{equation}\label{eq:mainequalityNSW}
	-\int_{\T}\mu k_\varepsilon\p_t\la\J\nabla u^\varepsilon,u_t^\varepsilon\nabla \varphi\ra d\vx-  \int_{\T}\frac{\p}{\p t}J(u^\varepsilon)\varphi d\vx=\int_{\T}\la\Hess(\varphi)\nabla u^\varepsilon,\J\nabla u^\varepsilon\ra d\vx.
\end{equation}

\section{Convergence and bounds of quantities related to the solution}\label{sec:NSW}

\subsection{Convergences of the Jacobian of the solution}

\begin{lemma}\label{thm: existence of vortices}
If the initial data $u^\varepsilon_0$ of \eqref{eq:NSW} satisfies Conditions \eqref{con:convergence of J varphi} and \eqref{Con:strong bd of energy}, then there exist $C^{1,1}$-paths $\vb_j:[0,T)\to\T$, such that
\begin{equation}\label{eq:converge of Jut and eut}
    J(u^\varepsilon(\vx,t))\to \pi\sum_{j=1}^{2N}d_j\delta_{\vb_j(t)},\quad k_\varepsilon e^\varepsilon(u^\varepsilon(\vx,t))\to\pi\sum_{j=1}^{2N}\delta_{\vb_j(t)}.
\end{equation}
\end{lemma}

\begin{proof}
We first denote
\begin{equation}\label{eq:def of r0}
    r_0=r(\va^0)=\frac{1}{4}\min_{j\ne k}|\va_j^0-\va_k^0|>0,
\end{equation}
and define
\begin{equation}\label{eq:define of Te}
T^\varepsilon:=\sup\left\{T>0:\|J(u_0^\varepsilon)- J(u^\varepsilon(t)) \|_{\M}\le \frac{\pi}{400}r_0\ \text{for any}\ t\in[0,T]\right\}.
\end{equation}
\eqref{Con:strong bd of energy}, \eqref{eq:hamiltonian conservation} and \eqref{eq:def of We} imply that there exists a constant $C>0$ such that for all $\varepsilon>0$,
\begin{equation}\label{con:weak bd of energy}
E^\varepsilon(u^\varepsilon(t))\le H_\mu^\varepsilon(u^\varepsilon(t))\le  2N \pi\log\left(\frac{1}{\varepsilon} \right)+C.
 \end{equation}
Then via the same method of the proof of (3.10)-(3.12) in \cite{zhu2023VDNW}, Theorem 1.4.4 in \cite{CollianderJerrard1999GLvorticesSE}, \eqref{eq:define of Te} and \eqref{con:weak bd of energy} imply that for each $0<t<T^\varepsilon,\varepsilon<\varepsilon_0$, there exist $\vb_j^\varepsilon(t)\in B_{r_0/2}(\va_j^0),j=1,\cdots,2N$,  such that
\begin{equation}\label{eq:distance bet Jut and dirac}
\begin{aligned}
\left\|J(u^\varepsilon(\vx,t))-\pi\sum_{j=1}^{2N}d_j\delta_{\vb_j^\varepsilon(t)} \right\|_{\M}=o(1),\\
\left\|k_\varepsilon e^\varepsilon(u^\varepsilon(\vx,t))-\pi\sum_{j=1}^{2N} \delta_{\vb_j^\varepsilon(t)}\right\|_{\M}=o(1),
\end{aligned}\end{equation}
\begin{equation}\label{eq:bound of ju, eu out balls}
\int_{\T_{3r_0/4}(\va^0)}e^\varepsilon(u^\varepsilon(\vx,t))d\vx \le C,\quad \|\vj(u^\varepsilon)\|_{L^1(\T)}\le C,\end{equation}
\begin{equation}\label{eq:lower bound of energy t}
E^\varepsilon(u^\varepsilon(t))\ge 2N\pi\log\frac{1}{\varepsilon}-C.
\end{equation}
\eqref{eq:lower bound of energy t} and \eqref{con:weak bd of energy} imply
\begin{equation}\label{eq:bound of ut}
    k_\varepsilon \|u^\varepsilon_t\|_{L^2(\T)}^2\le C
\end{equation}
and 
\begin{equation}\label{eq:distance bet eut and hut}
    \|k_\varepsilon e^\varepsilon(u^\varepsilon(\vx,t))-k_\varepsilon h^\varepsilon_\mu(u^\varepsilon(\vx,t)) \|_{\M}=\frac{\mu k_\varepsilon^2}{2}\left\||u_t^\varepsilon(\vx,t)|^2\right\|_{\M}\le Ck_\varepsilon.
\end{equation}

Then we estimate $|\vb_j^\varepsilon(t_2)-\vb_j^\varepsilon(t_1)|$: We can find $\eta\in C^\infty_0(B_{r_0}(\va_j^0))$ such that
\begin{equation*}
    \eta(\vx)=(\vx-\vb_j^\varepsilon(t_1))\cdot\frac{\vb_j^\varepsilon(t_2)-\vb_j^\varepsilon(t_1)}{|\vb_j^\varepsilon(t_2)-\vb_j^\varepsilon(t_1)|},\quad \vx\in B_{3r_0/4}(\va_j^0).
\end{equation*}
Combining \eqref{eq:distance bet Jut and dirac}, \eqref{eq:distance bet eut and hut}, \eqref{eq:derivative of h}, \eqref{con:weak bd of energy} and \eqref{eq:bound of ut}, we have
\begin{align*}
|\vb_j^\varepsilon(t_2)-\vb_j^\varepsilon(t_1)|=&\eta(\vb_j^\varepsilon(t_2))-\eta(\vb_j^\varepsilon(t_1))\\=&\int_{\T}k_\varepsilon\eta(\vx)(h^\varepsilon_\mu(u^\varepsilon(\vx,t_2))-h^\varepsilon_\mu(u^\varepsilon(\vx,t_1)))d\vx+o(1)\\
=&\int_{t_1}^{t_2}\int k_\varepsilon\eta \p_th^\varepsilon_\mu d\vx dt+o(1)\\=&-\int_{t_1}^{t_2}\int_{\T}k_\varepsilon\la\nabla \eta,\re(\overline{u_t^\varepsilon}\nabla u^\varepsilon)\ra d\vx dt+o(1)\\
\le&\|\nabla \eta\|_{L^\infty(\T)}k_\varepsilon\int_{t_1}^{t_2}\|\nabla u^\varepsilon\|_{L^2(\T)}\|u^\varepsilon_t\|_{L^2(\T)}dt+o(1)\le C|t_2-t_1|+o(1).
\end{align*}
Hence, there exist Lipschitz paths $\vb_j$'s and $T_0>0$ such that $\vb_j^\varepsilon(t)$ converges to $ \vb_j(t)$ in $[0,T_0]$ uniformly, which together with \eqref{eq:distance bet Jut and dirac}  implies \eqref{eq:converge of Jut and eut}. In particular, \eqref{con:convergence of J varphi} and \eqref{eq:converge of Jut and eut} imply
\begin{equation}\label{eq:bj0=aj0}
    \vb_j(0)=\va_j^0.
\end{equation}

Then, we prove that $\vb_j\in C^{1,1}([0,T_0],\T)$. 

For any unit vector $\bvec{\nu}\in \R^2$, $0<h<<1$ and $t\in[0,T_0],$ we can find $\varphi_1,\varphi_2\in C_0^\infty(B_{r_0}(\va_j^0) $ satisfying
\begin{equation*}
    \varphi_1=(\vx-\vb_j(t))\cdot\bvec{\nu},\quad \varphi_2=(\vx-\vb_j(t))\cdot\bvec{\nu}^\perp, \quad\vx\in B_{3r_0/4}(\va_j^0).
\end{equation*}
Then 
\begin{equation}\label{eq:JDf1=Df2}\nabla \J\varphi_1+\nabla\varphi_2=\bvec{0}\quad \text{in}\ B_{3r_0/4}(\va_j^0).\end{equation} \eqref{eq:mainequalityNSW} imply that
\begin{align}\label{eq:NSWBasicEuality}
    &-\int_{\R} \zeta^h(t-s)\int_{\T}\mu k_\varepsilon \p_t\la u_t^\varepsilon\nabla\varphi_1,\J\nabla u^\varepsilon\ra d\vx ds-  \int_{\R} \zeta^h(t-s)\int_{\T}\p_tJ(u^\varepsilon)\varphi_1d\vx ds\nn
    &\quad=\int_{\R} \zeta^h(t-s)\int_{\T}\la\Hess(\varphi_1)\nabla u^\varepsilon,\J\nabla u^\varepsilon \ra d\vx ds,
\end{align}
where $\zeta^h$ is defined by \cite{Jerrard1999GinzburgLandauwave}:
\begin{equation}\label{eq:def of zeta h}
    \zeta^h(t)=\begin{cases}
    0,&|t|\ge h,\\
    \frac{1}{h}-\frac{|t|}{h^2},&|t|\le h.
    \end{cases}
\end{equation}
By \eqref{eq:JDf1=Df2}, \eqref{eq:derivative of h}, \eqref{eq:bound of ju, eu out balls}, \eqref{eq:converge of Jut and eut}, \eqref{eq:bound of ut}  and \eqref{eq:distance bet eut and hut}, we obtain that the first term of the left hand-side of \eqref{eq:NSWBasicEuality} can be estimated as follows
\begin{align}\label{eq:first term of Basic eqn}
    &-\int_{\R} \zeta^h(t-s)\int_{\T}k_\varepsilon\p_t\la u_t^\varepsilon\nabla\varphi_1,\J\nabla u^\varepsilon\ra d\vx ds\nn 
    &\quad=-\int_{\R} \zeta^h(t-s)\int_{\T}k_\varepsilon\p_t\la u_t^\varepsilon\nabla\varphi_2,\nabla u^\varepsilon\ra d\vx ds\nn &\quad\quad+\int_{\R} \zeta^h(t-s)\int_{\T}k_\varepsilon\p_t\la u_t^\varepsilon(\J\nabla\varphi_1+\nabla \varphi_2),\nabla u^\varepsilon\ra d\vx ds\nn
    &\quad=\int_{\R} \zeta^h(t-s)\int_{\T}k_\varepsilon\p_t(\nabla\cdot(\re(\overline{u^\varepsilon_t}\nabla u^\varepsilon)))\varphi_2d\vx ds\nn&\qquad+O\left(\frac{1}{h}k_\varepsilon\|u_t^\varepsilon\|_{L^2(\T\times[t-h,t+h])}\|\nabla u^\varepsilon\|_{L^2(\T_{r_0}(\va^0)\times[t-h,t+h])}\right)\nn
    &\quad=\int_{\R}\zeta^h(t-s)\int_{\T}\p_t^2(k_\varepsilon h^\varepsilon_\mu(u^\varepsilon(s)))\varphi_2d\vx ds+O(\sqrt{k_\varepsilon})\nn
    &\quad=\int_{\T}k_\varepsilon\frac{1}{h^2}\left(h^\varepsilon_\mu(u^\varepsilon(t+h))-2h^\varepsilon_\mu(u^\varepsilon(t))+h^\varepsilon_\mu(u^\varepsilon(t-h)) \right)\varphi_2 d\vx+o(1)\nn
    &\quad=\pi\frac{\vb_j(t+h)-2\vb_j(t)+\vb_j(t-h)}{h^2}\cdot\bvec{\nu}^\perp+o(1).
\end{align}
Integrating by part with respect to $s$ and using \eqref{eq:distance bet Jut and dirac}, we obtain that the second term of the left hand-side of \eqref{eq:NSWBasicEuality} equals to
\begin{align}\label{eq:second term of Basic eqn}
      &\int_{\R}(\zeta^h)'(t-s)\int_{\T}J(u^\varepsilon(s))\varphi_1d\vx ds=\frac{ 1 }{h^2}\int_t^{t+h}\int_{\T}(J(u^\varepsilon(s))-J(u^\varepsilon(s-h)))\varphi_1d\vx ds\nn
    &\quad=\frac{d_j\pi  }{h^2}\int_t^{t+h}(\vb_j(s)-\vb_j(s-h))\cdot\bvec{\nu} ds+o(1)    \le C\|\vb\|_{C^1[t-h,t+h]}+o(1).
\end{align}
\eqref{eq:bound of ju, eu out balls} imply that the right hand-side of \eqref{eq:NSWBasicEuality} is less than
\begin{equation}\label{eq:RHS of basic}
    C\int_{\R}\zeta^h(t-s)\|\varphi_1\|_{C^2(\T)}\|\nabla u^\varepsilon\|^2_{\T_{3r_0/4}(\va^0)}ds\le C.
\end{equation}
Taking 
\begin{equation*}
    \bvec{\nu}^\perp=-\J\bvec{\nu}=\frac{\vb_j(t+h)-2\vb_j(t)+\vb_j(t-h)}{|\vb_j(t+h)-2\vb_j(t)+\vb_j(t-h)|},
\end{equation*}
and substituting \eqref{eq:first term of Basic eqn}, \eqref{eq:second term of Basic eqn} and \eqref{eq:RHS of basic} into \eqref{eq:NSWBasicEuality}, we have
\begin{equation*}
    \mu\pi\frac{|\vb_j(t+h)-2\vb_j(t)+\vb_j(t-h)|}{h^2}\le C\|\vb\|_{C^1[t-h,t+h]}+C+o(1),
\end{equation*}
which implies that $\vb\in C^{1,1}([0,T])$.

The above proof can be repeated as long as $\vb_j(T_0)\ne \vb_k(T_0)$ for any $1\le j<k\le 2N$. Hence, we can extend the existence interval of $\vb_j$, i.e. $[0,T_0]$, to $[0,T)$, with
\begin{equation*}
T:=\inf\left\{t|\, \vb_j(t)=\vb_k(t)\ \text{for some}\ 1\le j<k\le 2N \right\}.\qedhere
\end{equation*}
\end{proof}

\subsection{convergence of the current of the solution}

With $\vb=\vb(t)=(\vb_1(t),\cdots,\vb_{2N}(t))^T$ obtained in Lemma \ref{thm: existence of vortices}, we can define
\begin{equation}\label{eq:def of u*}
    u^*:=u^*(\vx,t)=H(\vx;\vb(t),\vq_*(\vb(t))),
\end{equation}
where $\vq_*(\vb(t))$ is continuous and satisfies
\begin{equation}\label{eq:def of qb}
    \vq_*(\vb(0))=\vq_0,\quad \vq_*(\vb(t))\in2\pi\sum_{j=1}^{2N}d_j\vb_j(t)+2\pi\Z^2.
\end{equation}

\begin{lemma}\label{lem:converge of j}
Assume $u^\varepsilon,\vb$ are the same as in Lemma \ref{thm: existence of vortices}. Then, for any $T_1<T$, $\rho=\min_{t\in[0,T_1]}r(\vb(t))$,
\begin{equation}\label{eq:converge of jut and j/||}
\vj(u^\varepsilon)\wto\vj(u^*)\ \text{in}\ L^1(\T\times[0,T_1]),\quad \frac{\vj(u^\varepsilon)}{|u^\varepsilon|}\wto \vj(u^*) \ \text{in}\ L^2(\T_\rho(\vb(t))\times[0,T_1]),
\end{equation}
with $u^*$ defined by \eqref{eq:def of u*}.
\end{lemma}
\begin{proof}

The proof is essentially the same as the proof of Lemma 3.2 in \cite{zhu2023VDNW}. The only difference is the verification of 
\begin{equation}
\nabla \cdot(\vj_*-\vj(u^*))=0,
\end{equation}
where $\vj_*$ satisfies 
\begin{equation}\label{eq:converge of j}
\vj(u^\varepsilon)\wto \vj_*\quad \text{in}\ L^1(\T\times[0,T_1]),\quad \vj(u^\varepsilon)=|u^\varepsilon|\frac{\vj(u^\varepsilon)}{|u^\varepsilon|}\wto \vj_*\quad \text{in}\ L^1(\T_\rho(\vb(t))\times[0,T_1]).
\end{equation}
The the proof of existence of $\vj_*$ is the same as in the proof of Lemma 3.2 in \cite{zhu2023VDNW}, so we omit it for brevity.

Theorem 1.4.1 in \cite{CollianderJerrard1999GLvorticesSE}, \eqref{con:weak bd of energy}, \eqref{eq:converge of Jut and eut} and \eqref{eq:def of e j J} imply
\begin{equation}\label{eq:bound of j/|| out balls}
\int_{0}^{T_1}\left\|\frac{\vj(u^\varepsilon)}{|u^\varepsilon|} \right\|_{L^2(\T_\rho(\vb(t)))}^2dt\le \int_{0}^{T_1}\int_{\T_\rho(\vb(t))}e^\varepsilon(u^\varepsilon(t))d\vx dt\le CT_1,
\end{equation}
\begin{equation}\label{eq:bound of L2norm}
\left\||u^\varepsilon|^2-1\right\|_{L^2(\T)}\le C\varepsilon\sqrt{|\log \varepsilon|}.
\end{equation}

For any $\varphi\in C^\infty_0(\T\times[0,T])$, \eqref{eq:converge of j}, \eqref{eq:div of ju}, \eqref{eq:bound of ut} and \eqref{eq:bound of L2norm} give
\begin{align*}
    \left|\int_{\T\times[0,T]}\nabla\varphi\cdot\vj_* d\vx dt\right|&=\lim_{\varepsilon\to0}\left|\int_{\T\times[0,T]}\nabla\varphi\cdot\vj(u^\varepsilon) d\vx dt\right|=\lim_{\varepsilon\to0}\left|\int_{\T\times[0,T]}\varphi\nabla\cdot\vj(u^\varepsilon) d\vx dt\right|\\
    &=\lim_{\varepsilon\to 0}\left|\int_{\T\times[0,T]} \varphi\frac{\p}{\p t}\left[\mu k_\varepsilon \im(\overline{u^\varepsilon}u^\varepsilon_t)-\frac{1}{2}  (|u^\varepsilon|^2-1)\right] d\vx dt\right|\\
    &=\lim_{\varepsilon\to 0}\left|\int_{\T\times[0,T]} \frac{\p}{\p t}\varphi\left[\mu k_\varepsilon \im(\overline{u^\varepsilon}u^\varepsilon_t)-\frac{1}{2}  (|u^\varepsilon|^2-1)\right] d\vx dt\right|\\
    &\le\lim_{\varepsilon\to 0} C(k_\varepsilon\|u_t^\varepsilon\|_{L^2(\T\times[0,T])}+\||u^\varepsilon|^2-1\|_{L^2(\T\times[0,T])}) =0,
\end{align*}
which imply that $\nabla \cdot\vj_*=0$. Then \eqref{eq:div of jH} and \eqref{eq:def of u*} imply
\begin{equation*}
    \nabla \cdot (\vj_*-\vj(u^*))=0.
\end{equation*}
Similarly, we have $\nabla\cdot(\J(\vj_*-\vj(u^*)))=0$. Then the proof can be completed by following the proof of Lemma 3.2 in \cite{zhu2023VDNW}.
\end{proof}

\subsection[Lower bounds of E(u) and the L2 norm of ut]{Lower bounds of $E^\varepsilon(u^\varepsilon(t))$ and the $L^2$-norm of $u^\varepsilon_t$ }
Lemma 2.3 in \cite{zhu2023quantized}, \eqref{eq:converge of Jut and eut} and \eqref{eq:converge of jut and j/||} imply that for any $0\le t_1<t_2<T$, and $\rho<\min_{t\in[t_1,t_2]}r(\vb(t))$, we have
\begin{equation}\label{eq:distance between ju and ju*}
0\le \limsup_{\varepsilon\to 0}\int_{t_1}^{t_2}\int_{\T_\rho(\vb(t))}\left(e^\varepsilon(|u^\varepsilon|)+\frac{1}{2}\left|\frac{\vj(u^\varepsilon)}{|u^\varepsilon|}-\vj(u^*) \right| \right)d\vx dt\le C\int_{t_1}^{t_2}\Sigma(t)dt,
\end{equation}
with
\begin{equation*}
    \Sigma(t):=\limsup_{\varepsilon\to 0}\left(E^\varepsilon(u^\varepsilon(t))-W_\varepsilon(\vb(t);\vq_*(\vb(t))) \right).
\end{equation*}
A direct corollary of \eqref{eq:distance between ju and ju*} is  $\Sigma(t)\ge 0$, i.e.
\begin{equation}\label{eq:refined lower bound of e}
\limsup_{\varepsilon\to 0}(E^\varepsilon(u^\varepsilon(t))- W_\varepsilon(\vb(t);\vq_*(\vb(t))))\ge 0.
\end{equation}

Corollary 7 in \cite{SandierSerfaty2004EstimateGinzburgLandau}, \eqref{eq:converge of Jut and eut}, \eqref{eq:bound of ut} and \eqref{con:weak bd of energy} imply 
\begin{equation}\label{eq:lower bound of L2 ut}
    \liminf_{\varepsilon\to 0}k_\varepsilon\int_{\T\times[t_1,t_2]}|u^\varepsilon_t(\vx,s)|^2d\vx ds\ge \pi\sum_{j=1}^{2N}\int_{t_1}^{t_2}|\dot{\vb}_j(s)|^2ds.
\end{equation}

\section{Proof of Theorem \ref{thm:dynamics NSW}}
\begin{proof}
Recall that $\va$ is the solution of \eqref{eq:NSWODE}, and  that $\vb$ is obtained in Lemma \ref{thm: existence of vortices}. Then we define
\begin{equation}\label{eq:def of zeta}
    \zeta(t)=\sum_{j=1}^{2N}\left(|\vb_j(t)-\va_j(t)|+\mu |\dot{\vb}_j(t)-\dot{\va}_j(t)|\right).
\end{equation}
And we can find $T_*<T$ such that 
\begin{equation}\label{eq:def of r2 and T2}
  \sup_{t\in[0,T_*]}\zeta(t)<r_*:=\inf_{t\in[0,T_*]}\min\{r(\va(t)),r(\vb(t))\}.
\end{equation}
For simplicity, we will use the following notations in this subsection: 
\begin{equation}\label{eq:simple notation}
\begin{array}{ll}
    W(\va(t))=W(\va(t);\vq_*(\va(t))),&\quad W(\vb(t))=W(\vb(t);\vq_*(\vb(t))),\\
    W_\varepsilon(\va(t))=W_\varepsilon(\va(t);\vq_*(\va(t))),&\quad W_\varepsilon(\vb(t))=W_\varepsilon(\vb(t);\vq_*(\vb(t))).
\end{array}
\end{equation}

We first verify that $\dot{\vb}(0)=\bvec{0}$. 
\eqref{eq:lower bound of L2 ut},  \eqref{eq:hamiltonian conservation}, \eqref{Con:strong bd of energy}, \eqref{eq:refined lower bound of e} and \eqref{eq:bj0=aj0} imply
\begin{align}\label{eq:vbdot=0}
\frac{\pi}{2h}\int_0^h|\dot{\vb}(t)|^2dt\le&\liminf_{\varepsilon\to0}\frac{1}{h}\int_0^h\int_{\T}\frac{k_\varepsilon}{2}|u_t^\varepsilon|^2 d\vx dt\nn=&\liminf_{\varepsilon\to0}\frac{1}{\mu h}\int_0^h\int_{\T}(h^\varepsilon(u^\varepsilon(t))-e^\varepsilon(u^\varepsilon(t)))d\vx dt\nn
\le&\frac{1}{\mu h}\int_0^h\left(W_\varepsilon(\va^0;\vq_0)-W_\varepsilon(\vb(t)) \right)dt\le C\frac{1}{\mu h}\int_0^h|\va_j^0-\vb_j(t)|dt\le Ch.
\end{align}
Letting $h\to 0$ on the both sides of \eqref{eq:vbdot=0}, we obtain $\dot{\vb}(0)=\bvec{0}$, which together with \eqref{eq:def of zeta}, \eqref{eq:initial ODE} and \eqref{eq:bj0=aj0} implies
\begin{equation}\label{eq:zeta0=0}
    \zeta(0)=0.
\end{equation}

For any $t<T_*$,
\begin{align}\label{eq:eq:first estimate of zeta'}
    \dot{\zeta}(t)\le&C\sum_{j=1}^{2N}|\dot{\vb}_j(t)-\dot{\va}_j(t)|+\sum_{j=1}^{2N}\left|\mu \ddot{\vb}_j(t)+  \J\dot{\vb}_j(t)-\mu \ddot{\va}_j(t)- \J\dot{\va}_j(t)\right|\nn
    \le&C\zeta+\frac{1}{\pi}\sum_{j=1}^{2N}\left|\nabla_{\vb_j}W(\vb(t))+\mu \pi \ddot{\vb}_j(t)+\pi   \J\dot{\vb}_j(t) \right|+\frac{1}{\pi}\sum_{j=1}^{2N}\left|\nabla_{\va_j}W(\va(t))-\nabla_{\vb_j}W(\vb(t)) \right|\nn
    \le&C\zeta+\frac{1}{\pi}|\bvec{A}_j(t)|,
\end{align}
where
\begin{equation*}
	\bvec{A}_j(t)=\nabla_{\vb_j}W(\vb(t))+\mu \pi \ddot{\vb}_j(t)+\pi   \J\dot{\vb}_j(t).
\end{equation*}

We may find $\varphi_1,\varphi_2\in C_0^\infty(B_{r_*}(\vb_j(t))) $ satisfying
\begin{equation*}
    \varphi_1=(\vx-\vb_j(t))\cdot\bvec{\nu},\quad \varphi_2=(\vx-\vb_j(t))\cdot\bvec{\nu}^\perp, \quad\vx\in B_{3r_*/4}(\vb_j(t)),
\end{equation*}
with 
\begin{equation}\label{eq:def of nuperp}
\bvec{\nu}^\perp=-\J\bvec{\nu}=\frac{\bvec{A}_j}{|\bvec{A}_j|}.
\end{equation}
Combining \eqref{eq:first term of Basic eqn}, \eqref{eq:second term of Basic eqn} and \eqref{eq:def of zeta h}, we have
\begin{align}\label{eq:aaaaa}
&\mu\pi\frac{\vb_j(t+h)-2\vb_j(t)-\vb_j(t-h)}{h^2}\cdot\bvec{\nu}^\perp+\frac{d_j\pi  }{h^2}\int_t^{t+h}(\vb_j(s)-\vb_j(s-h))\cdot\bvec{\nu}ds+o(1)\nn&\quad=-\int_{\R} \zeta^h(t-s)\int_{\T}\mu k_\varepsilon \p_t\la u_t^\varepsilon\nabla\varphi_1,\J\nabla u^\varepsilon\ra d\vx ds-  \int_{\R} \zeta^h(t-s)\int_{\T}\p_tJ(u^\varepsilon)\varphi_1d\vx ds.
\end{align}
We can derive from (4.17) in \cite{zhu2023quantized}, \eqref{eq:def of nuperp}, \eqref{eq:aaaaa}, \eqref{eq:mainequalityNSW}, \eqref{eq:production of jH and eta GP}, \eqref{eq:def of u*} that 
\begin{align}\label{eq:Adotnu}
    |\bvec{A}_j(t)|=&\bvec{A}_j(t)\cdot\bvec{\nu}^\perp\nn\quad=&\lim_{h\to 0}\left[\mu\pi\frac{\vb_j(t+h)-2\vb_j(t)+\vb_j(t-h)}{h^2}\cdot\bvec{\nu}^\perp\right.+\frac{d_j\pi  }{h^2}\int_t^{t+h}(\vb_j(s)-\vb_j(s-h))\cdot\bvec{\nu}ds\nn 
    &\qquad\left.-\int_{\R} \zeta^h(t-s)\int_{\T}\la\Hess(\varphi_1)\vj(u^*),\J\vj(u^*) \ra d\vx ds \right]\nn 
    \quad=&\lim_{h\to 0}\lim_{\varepsilon\to 0}\int_{\R} \zeta^h(t-s)\int_{\T}(\la\Hess(\varphi_1)\nabla u^\varepsilon,\J\nabla u^\varepsilon \ra-\la\Hess(\varphi_1)\vj(u^*),\J\vj(u^*) \ra)d\vx ds\nn
    \quad=&I_1+I_2,
\end{align}
with 
\begin{align*}
I_1=&\lim_{h\to 0}\lim_{\varepsilon\to 0}\int_{\R} \zeta^h(t-s)\int_{\T}(\la\Hess(\varphi_1)\nabla |u^\varepsilon|,\nabla |u^\varepsilon| \ra)d\vx ds\nn
&+\lim_{h\to 0}\lim_{\varepsilon\to 0}\int_{\R} \zeta^h(t-s)\int_{\T}\left(\la\Hess(\varphi_1)\left(\frac{\vj(u^\varepsilon)}{|u^\varepsilon|}-\vj(u^*) \right) ,\J\left(\frac{\vj(u^\varepsilon)}{|u^\varepsilon|}-\vj(u^*) \right) \ra\right)d\vx ds,\nn
I_2=&\lim_{h\to 0}\lim_{\varepsilon\to 0}\int_{\R} \zeta^h(t-s)\int_{\T}\left(\la\Hess(\varphi_1)\left(\frac{\vj(u^\varepsilon)}{|u^\varepsilon|}-\vj(u^*) \right) ,\J\vj(u^*) \ra\right)d\vx ds\nn
&+\lim_{h\to 0}\lim_{\varepsilon\to 0}\int_{\R} \zeta^h(t-s)\int_{\T}\left(\la\Hess(\varphi_1)\vj(u^*),\J\left(\frac{\vj(u^\varepsilon)}{|u^\varepsilon|}-\vj(u^*) \right) \ra\right)d\vx ds.
\end{align*}
\eqref{eq:converge of jut and j/||} immediately implies 
\begin{equation}\label{eq:I2=0}
I_2=0.
\end{equation}

Taking the derivative of $W(\va(s))+\frac{\mu \pi}{2}|\dot{\va}(s)|^2$ with respect to $s$ and noting \eqref{eq:NSWODE}, we have
\begin{align*}
&\frac{d}{dt}\left(W(\va(s))+\frac{\mu \pi}{2}|\dot{\va}(s)|^2\right)=\sum_{j=1}^{2N}\left(\dot{\va}_j(s)\cdot\nabla_{\va_j}W(\va(s))+\mu \pi \dot{\va}_j(s)\cdot\ddot{\va}_j(s)\right)\\
&\quad=\sum_{j=1}^{2N}\left(\dot{\va}_j(s)\cdot\nabla_{\va_j}W(\va(s))-\dot{\va}_j(s)\cdot\nabla_{\va_j}W(\va(s))-\pi d_j\dot{\va}_j(s)\cdot(\J\dot{\va}_j(s))\right)=0,
\end{align*}
which immediately implies that for any $s>0$, we have
\begin{equation}\label{eq:conserve W}
W(\va(s))+\frac{\mu \pi}{2}|\dot{\va}(s)|^2\equiv W(\va(0)).
\end{equation}
Then, combining \eqref{eq:hamiltonian conservation}, \eqref{Con:strong bd of energy}, \eqref{eq:def of We}, \eqref{eq:Lip W}  and \eqref{eq:conserve W}, we have
\begin{align}\label{eq:refined energy bound}
    E^\varepsilon(u^\varepsilon(s))=&H_\mu^\varepsilon(u^\varepsilon(s))-\int_{\T}\frac{\mu k_\varepsilon}{2}|u^\varepsilon_t|^2d\vx =E^\varepsilon(u^\varepsilon_0)+\frac{\mu k_\varepsilon}{2}\int_{\T}|u_1^\varepsilon|^2d\vx-\int_{\T}\frac{\mu k_\varepsilon}{2}|u^\varepsilon_t|^2d\vx\nn
    \le&W_\varepsilon(\va^0;\vq_0)-\int_{\T}\frac{\mu k_\varepsilon}{2}|u^\varepsilon_t|^2d\vx+o(1)\nn=&W_\varepsilon(\va(s))+\frac{\mu \pi}{2}|\dot{\va}(s)|^2-\int_{\T}\frac{\mu k_\varepsilon}{2}|u^\varepsilon_t|^2d\vx+o(1) \nn
    \le&W_\varepsilon(\vb(s))+C\zeta(s)+\frac{\mu \pi}{2}|\dot{\vb}(s)|^2-\int_{\T}\frac{\mu k_\varepsilon}{2}|u^\varepsilon_t|^2d\vx+o(1).
\end{align}
Combining \eqref{eq:refined energy bound}, \eqref{eq:distance between ju and ju*} and \eqref{eq:lower bound of L2 ut}, we obtain
\begin{align}\label{eq:estimate of I_1}
I_1\le &\lim_{h\to 0}\lim_{\varepsilon\to 0}\frac{C}{h}\int_{t-h}^{t+h}\int_{\T_{3r_*/4}(\vb(t))}\left|\frac{\vj(u^\varepsilon)}{|u^\varepsilon|}-\vj(u^*) \right|^2d\vx ds\nn
\le& \lim_{h\to 0}\lim_{\varepsilon\to 0}\frac{C}{h}\int_{t-h}^{t+h}\left(\zeta(s)+\frac{\mu \pi}{2}|\dot{\vb}(s)|^2-\int_{\T}\frac{\mu k_\varepsilon}{2}|u^\varepsilon_t(s)|^2d\vx \right)ds\le C\zeta(t).
\end{align}
Substituting \eqref{eq:estimate of I_1}, \eqref{eq:I2=0} and \eqref{eq:Adotnu} into \eqref{eq:eq:first estimate of zeta'}, we have
\begin{equation*}
	\dot{\zeta}(t)\le C\zeta(t).
\end{equation*}
Since $\zeta(0)=0$, we have $\zeta(t)=0$ for $t\in[0,T_*]$, which implies $\va(t)=\vb(t)$. Hence, the proof is completed.
\end{proof}

\section{Numerical results on the vortex dynamics}\label{sec:numerical}

The convergence of NLSW \eqref{eq:NSW} to NLS \eqref{eq:SE} as $\mu\to 0$ was widely studied numerically and analytically \cite{BaoCai2014NSW,BergeColin1995NSWinplasma,MachiharaNakanishiOzawa2002LimitOfNKG,Schoene1979NonrelativisticLimitKGandD}. And the convergence of the reduced dynamical law of NLSW \eqref{eq:NSW} on a bounded domain to the reduced dynamical law of NLS \eqref{eq:SE} when $\mu\to 0$ was established in \cite{Yu2011VortexDynamicsKleinGordon}. Based on Theorem \ref{thm:dynamics NSW} and results in \cite{zhubaojian2023quantized}, we study the convergence of the reduced dynamical law of NLSW \eqref{eq:NSW}  to the reduced dynamical law of NLS \eqref{eq:SE} as $\mu\to 0$ via numerical simulation. 

\begin{figure}[htp!]
\begin{center}
\begin{tabular}{ccc}
\includegraphics[height=4.3cm]{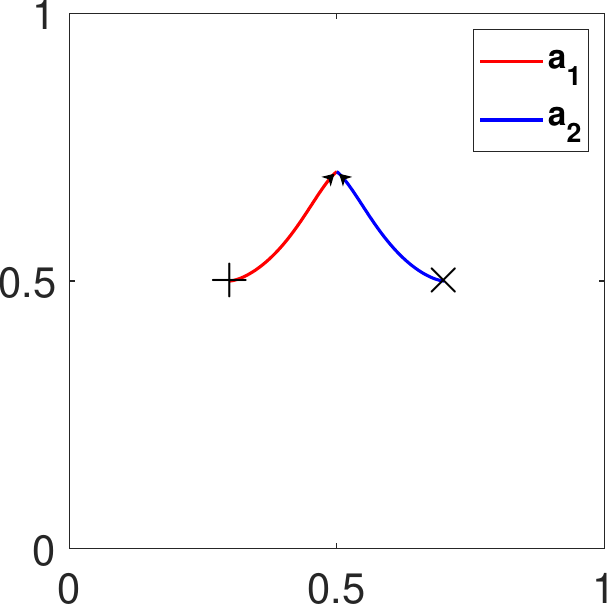}&&\includegraphics[height=4.3cm]{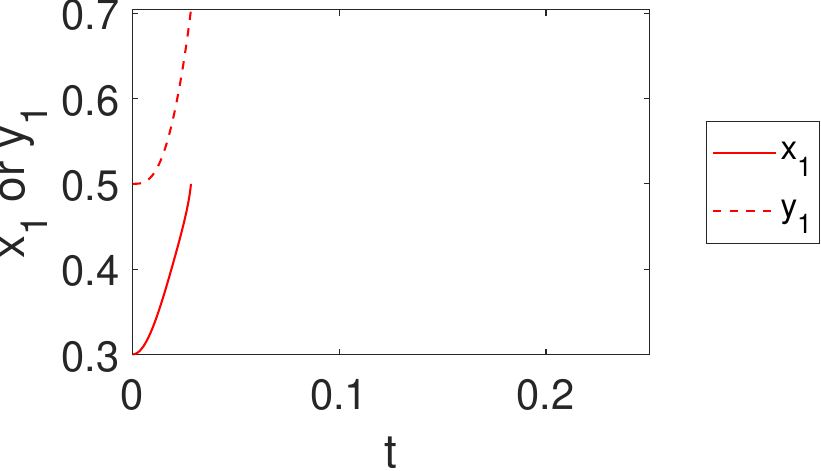}\\
\includegraphics[height=4.3cm]{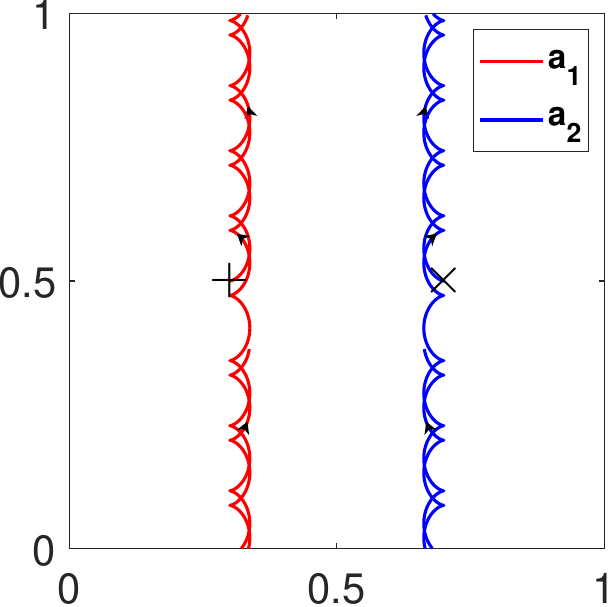}&&\includegraphics[height=4.3cm]{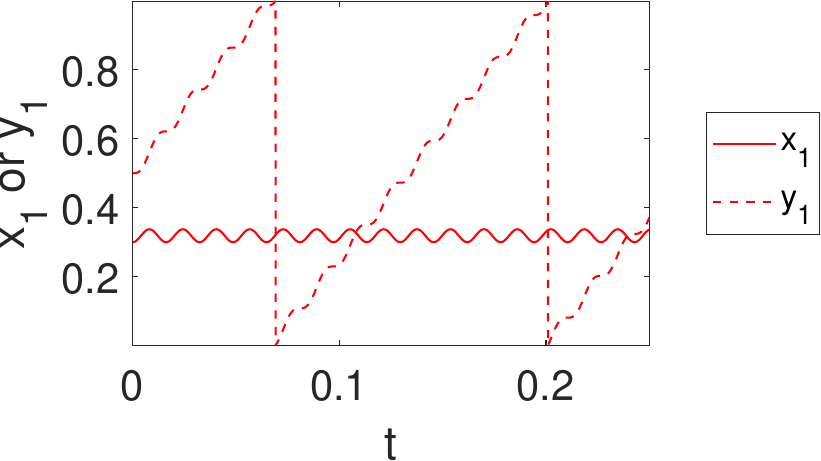}\\
\includegraphics[height=4.3cm]{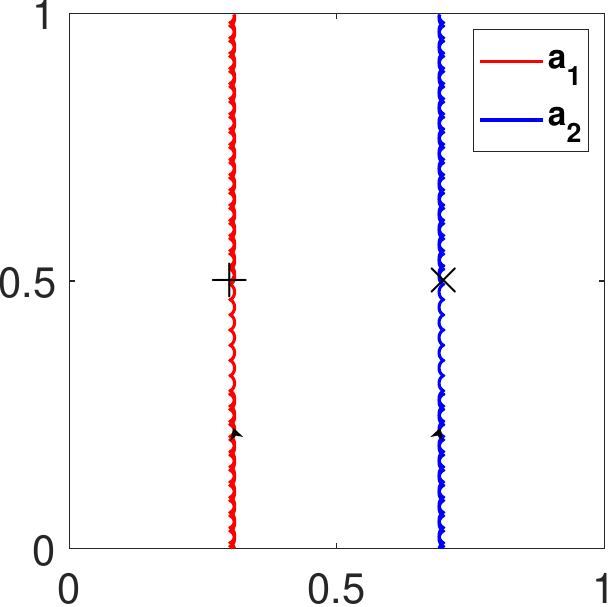}&&\includegraphics[height=4.3cm]{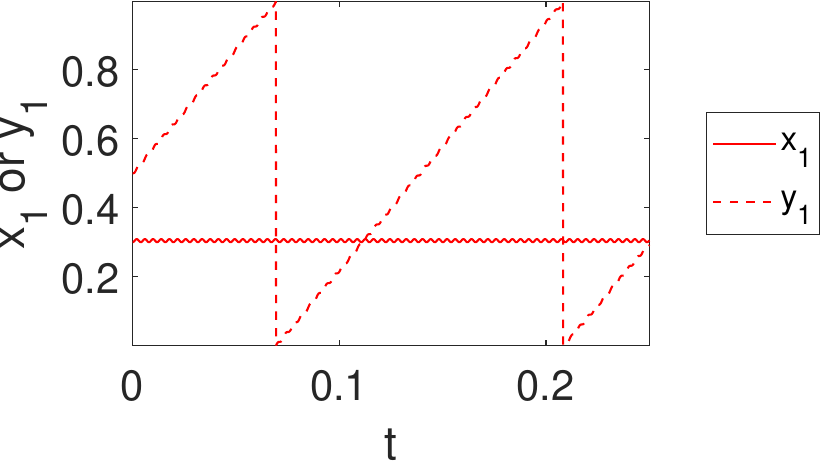}\\
\includegraphics[height=4.3cm]{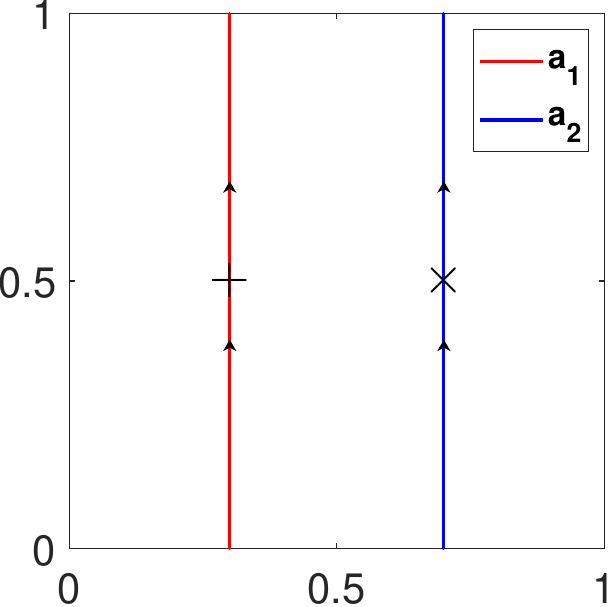}&&\includegraphics[height=4.3cm]{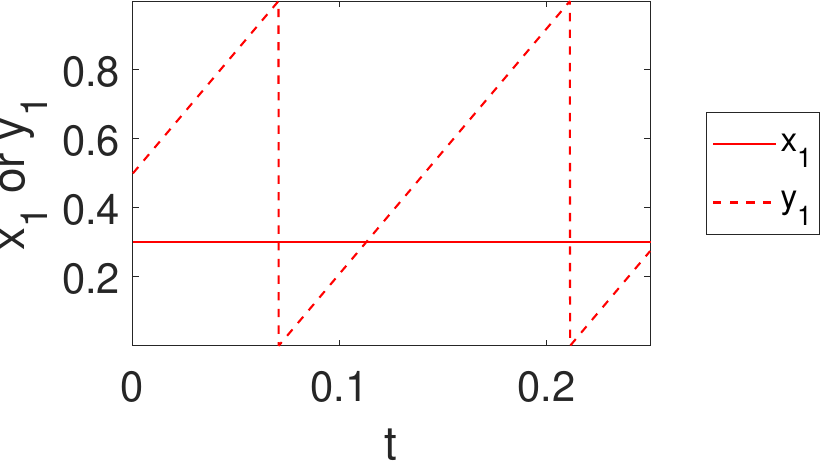}
\end{tabular}
\end{center}
\caption{Trajectories (left) and  values of $\va_1=(x_1,y_1)^T$ (right) of solutions of  \eqref{eq:NSWODE}  with different $\mu$ (from top to bottom line): (i) $\mu=1/100$; (ii) $\mu=1/400$; (iii) $\mu=1/1600$; (iv) $\mu=0$, i.e. the reduced dynamical law of NLS \ref{eq:SE}. The initial data are given by \eqref{eq:initial ODE} with $\va_1^0=(0.3,0.5)^T,\va_2^0=(0.7,0.5)^T,\vq_0=2\pi(\va_1^0-\va_2^0)$. Here , we use $+$ and $\times$ in the pictures to denote vortices with winding numbers $+1$ and $-1$, respectively.
}\label{fig:path}
\end{figure}

We focus on the case of the vortex dipole, i.e. $N=1$. We solve the ODE \eqref{eq:NSWODE} via the fourth-order Runge-Kutta method with time step $\Delta t=10^{-4}$. The results are shown in Figure \ref{fig:path}.

\section{Conclusion}

We have established the reduced dynamical law for quantized vortex dynamics of the nonlinear \SE{} with wave operator (NLSW) on the torus \eqref{eq:NSW} as the core size of vortex $\varepsilon\to 0$. When $\mu$ is a fixed positive number, the vortex motion law of NLSW is a mixed state of the vortex motions of the nonlinear wave equation and the nonlinear \SE{}. Finally, we investigate the convergence of the reduced dynamical law \eqref{eq:NSWODE} to the vortex motion law of the nonlinear \SE{} as $\mu\to 0$.

\section*{Acknowledgments}
This work was partially supported by the China Scholarship Council (Y. Zhu) and the National Natural Science Foundation of China (grant 12141103). Part of the work was done when the author was visiting National University of Singapore during 2021-2023 and the Institute for Mathematical Science in 2023. The author would like to express his sincere gratitude to Prof. Huaiyu Jian in Tsinghua University and Prof. Weizhu Bao in National University of Singapore for their guidance and encouragement.


\end{document}